\documentclass{IEEEtran}
\ifCLASSINFOpdf
\usepackage[pdftex]{graphicx}
\else
\usepackage[dvips]{graphicx}
\fi
\usepackage[cmex10]{amsmath}
\usepackage{cite}

\usepackage{setspace}
\usepackage{url}
\usepackage[colorlinks=true,linkcolor=black,anchorcolor=black,citecolor=black,urlcolor=black]{hyperref}
\usepackage{lipsum}
\usepackage{color}
\usepackage{array}
\usepackage{multirow}
\usepackage{makecell}
\usepackage{float}
\usepackage{setspace}
\usepackage{array}
\usepackage{supertabular}
\usepackage{supertabular,booktabs}
\usepackage{amssymb}
\usepackage{nomencl}
\usepackage{algorithm}
\usepackage{algpseudocode}
\allowdisplaybreaks[4]

\begin{document}
	\title{\huge Load Balancing in Low-Voltage Distribution Networks via Optimizing Residential Phase Connections}  
	\author{\small
		\IEEEauthorblockN{Bin Liu$^{1,2}$, Frederik Geth$^{1}$, Nariman Mahdavi$^{1}$, Jiangxia Zhong$^{3,4}$}\\
		\IEEEauthorblockA{$^{1}$Commonwealth Scientific and Industrial Research Organisation (CSIRO), Newcastle 2304, Australia\\
		$^{2}$School of EE\&T, The University of New South Wales, Sydney 2052, Australia\\
		$^{3}$School of Engineering, Royal Melbourne Institute of Technology, Melbourne 3000, Australia\\
		$^{4}$United Energy, Mount Waverley 3149, Australia\\
		\{brian.liu, frederik.geth, nariman.mahdavimazdeh\}@csiro.au, jiangxia.zhong@ue.com.au}} 
	\maketitle
	\newcommand\blfootnote[1]{%
		\begingroup
		\renewcommand\thefootnote{}\footnote{#1}%
		\addtocounter{footnote}{-1}%
		\endgroup}
	
	\begin{abstract}
        Unbalance issues in low-voltage distribution networks (LVDN) can be worsened by increasing penetration of residential PV generation if unevenly distributed among three phases. To address this issue, the phase-switching device (PSD) provides a viable and efficient method by dynamically switching customers to other phases. This paper further investigates how to optimize residential phase connections by controlling PSDs efficiently. The optimization problem is formulated as a mixed-integer non-convex programming (MINCP) problem considering relevant operational requirements of an LVDN based on the exact formulation of unbalanced three-phase optimal power flow (UTOPF). Unlike most heuristic algorithms and the  linearization techniques in our previous work, this paper proposes to solve the MINCP problem via an iteration-based algorithm after exact reformulations and reasonable approximations of some constraints. The proposed method is tested in a real LVDN and compared with the approach of Zhao et al. based on the well-known linear UTOPF formulation. Case studies based on the European low-voltage test feeder demonstrate the proposed method's efficiency in mitigating the network unbalance while ensuring network security and flexibility to deal with more controllable resources.
	\end{abstract}
		

	\section{Introduction}\blfootnote{\it The paper has been submitted to and accepted by IEEE PES Innovative Smart Grid Technologies Conference Asia, 2021 (ISGT Asia 2021).}
    The low-voltage distribution network (LVDN) in Australia is experiencing an increasing amount of challenges caused by high penetration of rooftop solar generation. 
    Amongst those challenges, addressing the phase unbalance has become  urgent in networks with PV systems unevenly distributed among the phases, particularly in the daytime when high irradiance occurs  \cite{RN61}. 
	The unbalance issue can cause lots of operation problems, e.g. malfunction of the protection system, increased power loss, low power supply quality and inefficient usage of electrical equipment  \cite{RN103}.
	The challenge is very significant in Australia as it is experiencing remarkable renewable energy deployment \cite{RN192}. Power utilities in Australia usually run an extensive four-wire grid along the streetscape, and most residential customers are powered by single or three phases. 
	
	To effectively address unbalance issue in LVDN, as well as improving its operational feasibility as revealed in our previous research in \cite{BLiu-TUPF}, one potentially cost-efficient option is using the phase-switching device (PSD) to dynamically switch the phase connections of residential customers equipped with PSDs, where the PSD controller is installed at the secondary side of the distribution transformer (DT), and controllable switches are distributed in the network. As discussed in \cite{psd-knowledge}, although PSDs could have some negative impacts on some residential appliances, most of the appliances works normally during the laboratory tests due to the advanced electronic technologies to detect the zero-crossing point of current wave forms.
	
	The controllable switches communicate with the controller wirelessly to switch their phase connections after receiving the optimized strategies. Voltages and currents of customers are assumed to be measurable by smart meters and available for the PSD controller via communication with the operation centre \cite{Liu2021}. The voltage at the secondary side of DT is also assumed to be available for the PSD controller through fast measurements in the PSD controller \cite{psd-knowledge}.
	

	Given the available data, the critical question is how to optimally coordinate the PSDs to minimize the unbalance level in the network. Reported methods in this field are mainly based on heuristic algorithms supported by black-box unbalanced three-phase power flow (UTPF) programs \cite{RN71,RN98}. The disadvantage of this type of approach is that the optimality cannot be guaranteed and its inflexibility in the problem formulation when the investigated component is not available for the UTPF solver. Other linear approximation \cite{RN103} or mixed-integer linear programming (MILP)-based algorithms \cite{RN44,RN105} either treat all customers as constant current loads  or do not consider the network model.

	A comprehensive formulation of the problem, which is based on unbalanced three-phase optimal power flow (UTOPF) and considers all operational constraints in the decision-making process, along with the linearization-based approach to solve the challenging mixed-integer non-convex programming (MINCP) problem is presented in our previous work \cite{PRD_LVDN_01}. However, the method will needs to be extended to support controllable resources, e.g., PV inverters, demand response or electric vehicles (EVs), motivating us to develop a more flexible yet efficient enough algorithm to solve the problem. 
	
	Compared with our previous work, this paper focuses on optimal coordination of PSDs in real-time operation, while \cite{PRD_LVDN_03} investigates the optimal locations to place PRDs, \cite{PRD_LVDN_01} focuses on the day-ahead scheduling of both PSDs and static VAR generator based on predicted demand profiles, and \cite{PRD_LVDN_02} mainly discussed the optimal control PRDs when the network cannot be monitored all all nodes, which leads to an sensitivity-based approach. Moreover, an iteration-based algorithm is presented and compared with the complete linearization-based method in our previous work \cite{PRD_LVDN_01} and the method based on the famous UTOPF that is based on linear branch flow model \cite{RN64}. The rest of the paper is structured as follows. The problem formulation is presented in Part II, followed by the solution techniques in Part III. Case studies based on the European low-voltage test feeder (European LVDN) is performed in Part IV, and the paper is concluded in Part V.  
	
	\section{Problem Formulation}
	In the formulation, $\mathcal {C}_i$ will be used to represent the set of customers connected to node $i$, $\mathcal{F}_i$ and $\mathcal{X}_i$ as the set of customers with and without PSDs installed (denoted as fixed and adjustable customers), respectively. Obviously, there are $\mathcal{F}_i\cap \mathcal{X}_i=\emptyset$ and $\mathcal{F}_i\cup \mathcal{X}_i=\mathcal{C}_i$. Moreover, $V$ represent the nodal voltage in the main feeder and $X$ and $Y$ are the real and imaginary parts of $V$; $I$ represents the current with $J$ and $W$ being the real and imaginary parts, respectively. The subscripts  $\phi/\psi \in \{a,b,c\}$ 
	represent the phase labels, node indices are $i$ and $k$; lines $ik$ connect from $i$ to $k$; the customer index is $j$. Other parameters or variables will be explained right after their appearances.
	
	Moreover, we assume that the voltage phasor of the root node, i.e. the primary side of DT, is a known and available parameter. All residential demands are assumed to be with constant active and reactive powers, which, through the installed smart meters, are available for both the operation centre and the PSD controller. The probelm can be formulated as follows. 
	\begin{subequations}
	\footnotesize
		\begin{eqnarray}
		\label{obj-I01}
		F=\min{\pi+M_{b}[\sum_{i}{(\tau^{-}_{i}+\tau^{+}_{i})}+\sum_{\phi}{\rho_{\phi}}+\sum_{i}{\omega^-_{i}}]}\\
		\label{obj-I02}
		|P_{\phi,i_0k_0}-P_{\psi,i_0k_0}|\le\pi, |Q_{\phi,i_0k_0}-Q_{\psi,i_0k_0}|\le \pi~\forall \phi,\forall \psi\\
		\label{ol-1}
		V_{\phi,i}-V_{\phi,k}=\sum\nolimits_\psi{Z_{ik}^{\phi\psi}I_{\phi,ik}}~\forall \phi,\forall ik\\
		\label{kcl-1}
		\sum_{n:n\rightarrow i}{I_{\phi,ni}}-\sum_{k:i\rightarrow k}{I_{\phi,ik}}=\sum_{j\in \mathcal {C}_i}{I_{\phi,i,j}}~\forall \phi,\forall i\neq x\\
		\label{kcl-2}
		I_{\phi,i,j}=\varepsilon_{\phi,i,j}(P^n_{i,j}-\text{j}Q^n_{i,j})/V_{\phi,i}^H~\forall \phi,\forall i,\forall j\\
		\label{kcl-4}
		\varepsilon_{\phi,i,j}\in \{0,1\}~\forall \phi,\forall i,\forall j\\
		\label{kcl-5}
		\sum\nolimits_\phi{\varepsilon_{\phi,i,j}=1}~\forall i,\forall j\\
		\label{vlimit-1}
		V_{\phi,i}^\text{min}-\tau_i^-\le |V_{\phi,i}|\le V_{\phi,i}^\text{max}+\tau_i^+~\forall \phi,\forall i\\
		\label{vroot}
		V_{\phi,i_0}=V_{\phi}^0~\forall \phi\\
		\label{negV-1}
		3V^-_{i}=[1,\chi,\chi^2][V_{a,i},V_{b,i},V_{c,i}]^T~\forall i\\
		\label{negV-2}
		|V^-_{i}|\le \nu+\omega_i^-~\forall i\\
		\label{Ilimit}
		|S_{\phi,i_0k_0}|=|V_{\phi}^0||I_{\phi,i_0k_0}|\le |V_{\phi}^0|(I_{\phi,i_0k_0}^\text{max}+\rho_\phi)=S_{i_0k_0}^{\phi,\text{max}}~\forall \phi\\
		\label{slack}
		\tau_i^-\ge 0,\tau_i^+\ge 0,\rho_\phi\ge 0,\omega_i^-\ge 0
		\end{eqnarray}
	\end{subequations}
	where
	$P_{\phi,i_0k_0}$ and $Q_{\phi,i_0k_0}$ are active and reactive powers running through the DT, which can be expressed as $P_{\phi,i_0k_0}+\text{j}Q_{\phi,i_0k_0}=V_{\phi,i_0}I^H_{\phi,i_0k_0}$ with $V_{\phi,i_0}$ and $I_{\phi,i_0k_0}$ being the voltage of root node $i_0$ and the current of line $i_0k_0$ from node $i_0$ to node $k_0$ in phase $\phi$.
	Superscript $H$ indicates the Hermitian adjoint (conjugate transpose);
	$Z_{ik}^{\phi\psi}$ is the mutual impedance between phase $\phi$ and $\psi$ of line $ik$; Moreover, $Z_{ik}$ will be used to represent the impedance matrix of line $ik$.
	$\varepsilon_{\phi,i,j}=\alpha_{\phi,i,j}$ for adjustable customers and $\varepsilon_{\phi,i,j}=\mu_{\phi,i,j}$ for fixed customers with $\alpha_{\phi,i,j}$ and $\mu_{\phi,i,j}$ being introduced binary variables indicating whether the $j^\text{th}$ customer is connected to phase $\phi$ of node $i$ and the initial phase connection of the same customer, respectively.
	$P^n_{i,j}/Q^n_{i,j}$ is net active/reactive demand of customer $j$ at node $i$. 
	$V_{\phi,i}^\text{min}/V_{\phi,i}^\text{max}$ is the lower/upper voltage magnitude (VM) limit of $V_{\phi,i}$.
	$V_{\phi}^0$ is the known voltage of the root node at phase $\phi$.
	$I_{\phi,i_0k_0}^\text{max}$ is the upper current magnitude limit of $I_{\phi,i_0k_0}$.
	$\chi=e^{-\text{j}2\pi/3}$ and $\nu$ is the negative sequence limit (a metric for voltage unbalance), which is 1\% in this paper.
	Non-negative slack variables are also introduced and constrained as \eqref{slack} to ensure the feasibility of the problem.
	
	In the problem specification, the objective is to minimize the unbalance level of active/reactive power running through the DT, which is defined by \eqref{obj-I02}, plus the penalty of the sum of non-negative slack variables.
	The Ohm's law for each line and Kirchhoff's current law for each node are satisfied by \eqref{ol-1} and \eqref{kcl-1}-\eqref{kcl-5}, respectively.
	Secure operation requirements, including VM constraints, voltage information of root node, negative sequence constraints and DT capacity constraints are formulated as \eqref{vlimit-1}, \eqref{vroot}, \eqref{negV-1}-\eqref{negV-2} and \eqref{Ilimit}, respectively.

	Obviously, the problem, denoted as OPSD, belongs to MINCP due to introduced integer variables, and non-convex constraints in \eqref{kcl-2} and \eqref{vlimit-1}. Solution approach will be discussed in the next section.
	
	\section{Solution Approach}\label{soltech}
	To effectively solve the problem, the non-convex parts in the formulation need to be  addressed. Specifically, constraints can be divided into three categories: non-linear but convex (NLC) constraints including \eqref{obj-I02}, \eqref{vlimit-1}$_{ul}$, \eqref{Ilimit} and \eqref{negV-2}, general non-convex (GNC) constraints including \eqref{kcl-2},\eqref{vlimit-1}$_{ll}$ and all other constraints as linear (LIN) constraints\footnote{$ul$ and $ll$ represent the upper and lower limits, respectively.}. 
	
	
	\subsection{Previous Work Review}
	Generally, NLC constraints can be kept in the formulation or efficiently approximated by a set of linear constraints(e.g. \cite{linear-SOC}) and, together with LINs, can be addressed by existing commercial solvers. In this paper, NLCs are reformulated as linear or second-order cone constraints before dealing with GNC constraints. 
	
	For the GNC constraints, \eqref{vlimit-1}$_{ll}$ is linearly approximated based on the fact that voltage angle ranges of all nodes in each phase are sufficiently small as demonstrated in many references, e.g., \cite{RN38,RN66}. Specifically, assuming VM limit at node $i$ is $[V^\text{min}_{i},V^\text{max}_{i}]$ for any phase and VA limit in phase $\phi$ for any node is $[\delta^\text{min}_\phi,\delta^\text{max}_\phi]$, where $\delta^\text{min}_\phi=\delta_\phi-\Delta \delta,\delta^\text{max}_\phi=\delta_\phi+\Delta\delta$, lower limit in \eqref{vlimit-1} can be approximately linearized as 
	\begin{eqnarray}
	\footnotesize
	\label{vm-lower-1}
	\small
	    X_{\phi,i}\cos{\delta_\phi}+Y_{\phi,i}\sin{\delta_\phi}\ge V^\text{min}_{\phi,i}-\tau_i^-
	\end{eqnarray}
	
    Therefore, how to deal with GNC constraint \eqref{kcl-2} is critical in the efficiency of solving the MINCP problem. In our previous work, the constraint is completely linearized after approximating $1/V^H$ by an affine function of $X$ and $Y$, which leads to \eqref{pb-1}, and reformulating the bilinear terms via McCormick Envelopes \cite{mccormick}.
    \begin{eqnarray}
    \footnotesize
	\label{pb-1}
	    I_{\phi,i,j}=\varepsilon_{\phi,i,j}(P^n_{i,j}-\text{j}Q^n_{i,j})[k^X_{\phi}X_{\phi,i}+k^Y_{\phi}Y_{\phi,i}+b^X_{\phi}\nonumber\\
	    +\text{j}(h^X_{\phi}X_{\phi,i}+h^Y_{\phi}Y_{\phi,i}+b^Y_{\phi})]~\forall \phi, \forall i,\forall j
	\end{eqnarray}
	where $k^X_\phi,k^Y_\phi,b^X_\phi$ and $h^X_\phi,h^Y_\phi,b^Y_\phi$ are parameters to be fitted for phase $\phi$ before the optimization process based on historical or empirical data.
	
	More details of the above formulation, which includes \eqref{obj-I01}-\eqref{kcl-1}, \eqref{kcl-4}-\eqref{kcl-5}, \eqref{vlimit-1}$_{ul}$, \eqref{vroot}-\eqref{slack} and \eqref{vm-lower-1}-\eqref{pb-1}, and is denoted as linearized voltage method (LINV-M) and leads to a mixed-integer second-order cone programming (MISOCP) problem, can be found in \cite{PRD_LVDN_01}. LINV-M will be will be used for comparison purposes in this paper. Moreover, it is should be noted that LINV-M will not longer apply when controllable resources from residential customers are considered, which, in other words, implies $P^n_{i,m}$ and/or $Q^n_{i,j}$ will be variables, even the aforementioned solution techniques are used due to the introduction of extra  bilinear terms. 
	
	\subsection{Approach Based on Well-Known Linearized UTOPF}
	To address GNC constraints, another viable approach is based on the well-known linearized UTOPF formulation proposed in \cite{RN64}. The discussed UTOPF is based branch power flow model, neglects the power losses and assumes phase angles of voltages among three phases are exactly with 120$^\circ$ difference. Assuming $V_i=[V_{a,i},V_{b,i},V_{c,i}]^T$, $v_i=V_iV_i^H$, $S_{ik}=[S_{a,ik},S_{b,ik},S_{c,ik}]^T$, $S_{\phi,ik}=P_{\phi,ik}+\text{j}Q_{\phi,ik}$, $S^n_{i,j}=P^n_{i,j}+\text{j}Q^n_{i,j}$ and $\beta=[1,\chi^2,\chi;\chi,1,\chi^2;\chi^2,\chi,1]$, \eqref{ol-1}-\eqref{kcl-2}, \eqref{vlimit-1}-\eqref{negV-2} based on the linear UTOPF is reformulated as,
	\begin{subequations}
	\footnotesize
		\begin{eqnarray}
		\label{vdrop-1}
		v_i-v_k=\beta\text{diag}(S_{ik})Z_{ik}^H+Z_{ik}(\beta\text{diag}(S_{ik}))^H~\forall ik\\
		\label{kclP-1}
		\sum_{n:n\rightarrow i}{S_{\phi,ni}}-\sum_{k:i\rightarrow k}{S_{\phi,ik}}=\sum_{j\in \mathcal {C}_i}{S_{\phi,i,j}}~\forall \phi,\forall i\neq x\\
		\label{kclP-2}
		S_{\phi,i,j}=\sum\nolimits_{\phi}\varepsilon_{\phi,i,j}S^n_{i,j}~\forall \phi,\forall i,\forall j\\
		\label{vlimitP-1}
		(V_{\phi,i}^\text{min})^2-\tau_i^-\le \text{diag}(v_{i})_\phi\le (V_{\phi,i}^\text{max})^2+\tau_i^+~\forall \phi,\forall i\\
		\label{vrootP}
		v_{i_0}=V^0(V^0)^H\\
		\label{negVP-2}
		9v^-_{i}=[1,\chi,\chi^2]v_i[1,\chi,\chi^2]^T, v^-_{i}\le \nu^2+\omega_i^-~\forall i
		\end{eqnarray}
	\end{subequations}
	
	The above formulation, together with \eqref{obj-I01}-\eqref{obj-I02}, \eqref{kcl-4}-\eqref{kcl-5}, \eqref{vroot}, and \eqref{Ilimit}-\eqref{slack}, is a linear programming problem, which is denoted as the linear branch flow model method (LBFM-M). 
	
	\subsection{The Iteration-based Algorithm}
	Noting that if $V^H_{\phi,i}$ is fixed and known, the GNC constraint \eqref{kcl-2} becomes linear, and OPSD it therefore an efficient solvable MILP problem. This motivates us to solve the problem by iteratively updating the fixed nodal voltages in the formulation. 
	With fixed nodal voltages, say $V^0_{\phi,i}$, \eqref{kcl-2} can be expressed as the following linear constraints.
	\begin{eqnarray}
	\footnotesize
	\label{kcl-L-2}
	I_{\phi,i,j}=\varepsilon_{\phi,i,j}(P^n_{i,j}-\text{j}Q^n_{i,j})V^0_{\phi,i}/|V^0_{\phi,i}|^2~\forall \phi,\forall i,\forall j
	\end{eqnarray}
	
	With \eqref{kcl-L-2} and other NLC/LIN constraints, OPSD is also reformulated as a MISOCP problem, which includes \eqref{obj-I01}-\eqref{kcl-1}, \eqref{kcl-4}-\eqref{kcl-5}, \eqref{vlimit-1}$_{ul}$, \eqref{vroot}-\eqref{slack}, \eqref{vm-lower-1} and \eqref{kcl-L-2}, and is denoted as fixed voltage method (FIXV-M) throughout the context. The specific algorithm is summarised as Algorithm \ref{alg-1} and it should be noted that the algorithm may lead to different optimal solutions with various $K$. 
	\begin{algorithm}
	\footnotesize
		\caption{\small\textit{FIXV-M to solve OPSD}}
		\label{alg-1}
		\begin{algorithmic}[1]
			\State Initialise $\Delta V=100$, $V^{0}_{\phi,i}=V^\text{fix}_{\phi,i}~(\forall \phi,\forall i)$, $\epsilon_V=10^{-4}$, $k=1$ and $K$, which is the maximum allowable iterations.
			\While {$\Delta V>\epsilon_V$ \& $k\le K$}
			\State Solve FIXV-M with fixed $V^0_{\phi,i}~\forall \phi,\forall i$. Denote the optimum of $V_{\phi,i}$ as $V^\text{opt}_{\phi,i}$.
			\State Update $\Delta V=\max_{\phi,i}{(\big|V^\text{opt}_{\phi,i}-V^{0}_{\phi,i}\big|)}$ and $V^{0}_{\phi,i}=V^\text{opt}_{\phi,i}$.
			\State $k=k+1$.
			\State Report the optimal solution as $V^\text{opt}_{\phi,i}~(\forall \phi,\forall i)$.
			\EndWhile
		\end{algorithmic}
	\end{algorithm}
	
	Several remarks on the algorithm are given below.
	\begin{enumerate}
		\item The optimal strategy provided by FIXV-M is expected to be as accurate as LINV-M because no linearization is employed in addressing the power balance equations. In addition, both LINV-M and FIXV-M are expected to be more accurate than LBFM-M because LBFM-M takes stronger linearizations and assumptions. 
		The accuracy of the three methods will be further investigated via the UTPF algorithm presented in \cite{RN37} after the optimal PSD phase connections are reported.
		\item Regarding the computational efficiency, both LINV-M and FIXV-M are expected to be less efficient than LBFM-M because LBFM-M is an efficiently solvable MILP while both LINV-M and FIXV-M are formulated as less-efficient MISOCP. Moreover, the computational time of FIXV-M will also depends on the value of $K$. However, when PSD phase connections are dynamically updated in the network operation, the computational efficiency of FIXV-M can be improved via setting $K$ as 1, which implies the algorithm will be run only once regardless of $\Delta V$ and $V^0_{\phi,i}$  as values collected from the network monitoring system or from running UTPF before doing the optimization.
		\item Both FIXV-M and LBFM-M are more flexible than LINV-M, noting that LINV-M is only viable when residential demands are fixed. In other words, the method will no longer apply if $P^n_{i,j}$ and/or $Q^n_{i,j}$ are variables.
		\item When $K$ is a large enough number, the convergence of Algorithm \ref{alg-1} cannot be guaranteed theoretically and needs more effort in our future work. Meanwhile, the maximum number of iterations, i.e. $K$, is introduced to ensure a solution can be provided in a reasonable number of iterations. However, simulation results showed that this method is highly accurate, which will be discussed in detail in the next section.  
		\item All MILPs/MISOCPs are solved by Gurobi 9.1 \cite{gurobi} on a computer with Intel i7-8550U 1.8 GHz CPU, 16 GB memory, and the time limit for the solver is set to 10 minutes.
	\end{enumerate}
	

	\section{Case Study}
    \subsection{Case setup}
\begin{figure}[htb!]
		\centering\includegraphics[scale=0.20]{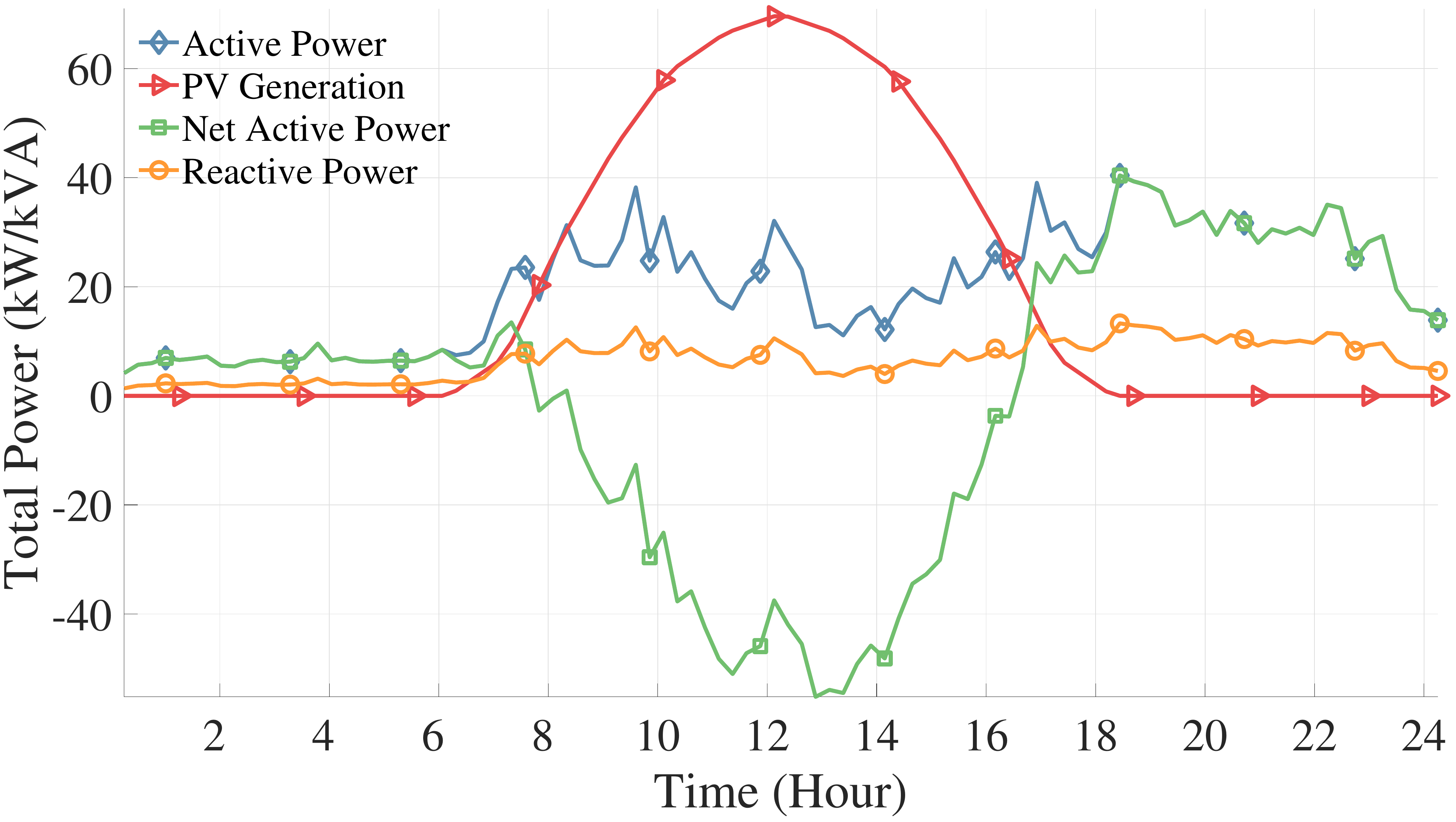}
		\caption{System loads in the real LVDN.}
		\label{fig-system-demand}
	\end{figure}

	The European LVDN, where the topology data and residential demand data can be found in \cite{EuroLVDN}, will be studied in this section. Among the 55 customers, we assume 10 of them (customer index: 5,9,15,18,20,26,30,37,45 and 50) are assigned PV systems, and 10 of them (customer index: 2,8,23,24,29,32,33,35,38 and 53) are assigned PSDs. Moreover, the installed capacity of each PV is assumed to be 7 kW, leading to a total capacity of 70 kW. The total system demands and PV generations with a time resolution of 15 minutes of a day are presented in Fig.\ref{fig-system-demand}, where the system exports power to the upper-level power grid around the midday due to high PV generation, and experiences peak demand around 18:30. Moreover, none of the PV customers have PSDs, to avoid unexpected shutdowns of PV inverter, as the impact of switching phases on the PV inverter operation needs more investigations \cite{psd-knowledge}.
	
	In the formulation, the voltage of the root node is a known parameter, which is set as 
	$$V^0=[1.05e^{j0},1.05e^{-j\frac{2\pi}{3}},1.05e^{j\frac{2\pi}{3}}]^T$$
	and $\nu^-$ is set as 1\% for all nodes. Typically, the $V^\text{min}$ and $V^\text{max}$ are set as 0.94 p.u. and 1.1 p.u., respectively. In the formulated problems, $M_b$ is set as 500 and the DT capacity is 200 KVA. Other parameters can be found in \cite{EuroLVDN}. 
	
	Two cases, which are outlined as follows, will be studied in this paper. 
	\begin{enumerate}
		\item The Base Case: All the three methods, i.e., LINV-M, LBFM-M and FIXV-M, will be tested. Moreover, two sub-cases of FIXV-M, which are described as follows, will be studied.
		\begin{itemize}
    		\item FIXV-MC (with cold-start point): $K$ is set as 3, and $V^0_{\phi,i}$ is set as $V^0_\phi$ for all periods throughout the day. In other words, Algorithm \ref{alg-1} will always be started from a cold-start point, i.e., setting $V^\text{fix}_{\phi,i}=V^0_{\phi}$, where $V^0_{\phi}$ is the voltage at the network's root node.
    		\item FIXV-MW (with warm-start point): $K$ is set as 1, and $V^0_{\phi,i}$ is set as the values from the network voltage data, which can be collected from the network monitoring system or via running UTPF, before doing the optimization.
    	\end{itemize}
		\item The Case with Controllable PVs: In this case, only LBFM-M and FIXV-MW will be studied and reactive powers of PV inverters are controllable.
	\end{enumerate}
	
	\subsection{The Base Case}
	For the Base Case, the solver averagely takes 24.72, 8.22, 9.18 and 3.32 seconds for each period under LINV-M, LBFM-M, FIXV-MC and FIXV-MW, respectively. As expected, LBFM-M and FIXV-MW are the most efficient on computational efficiency. It is interesting that although FIXV-MC solves a series of MISOCPs, it still takes less time than LINV-M, implying FIXV-MC would be a better choice than LINV-M if it brings comparable network performance improvements.
	
\begin{figure}[htb!]
		\centering\includegraphics[scale=0.2]{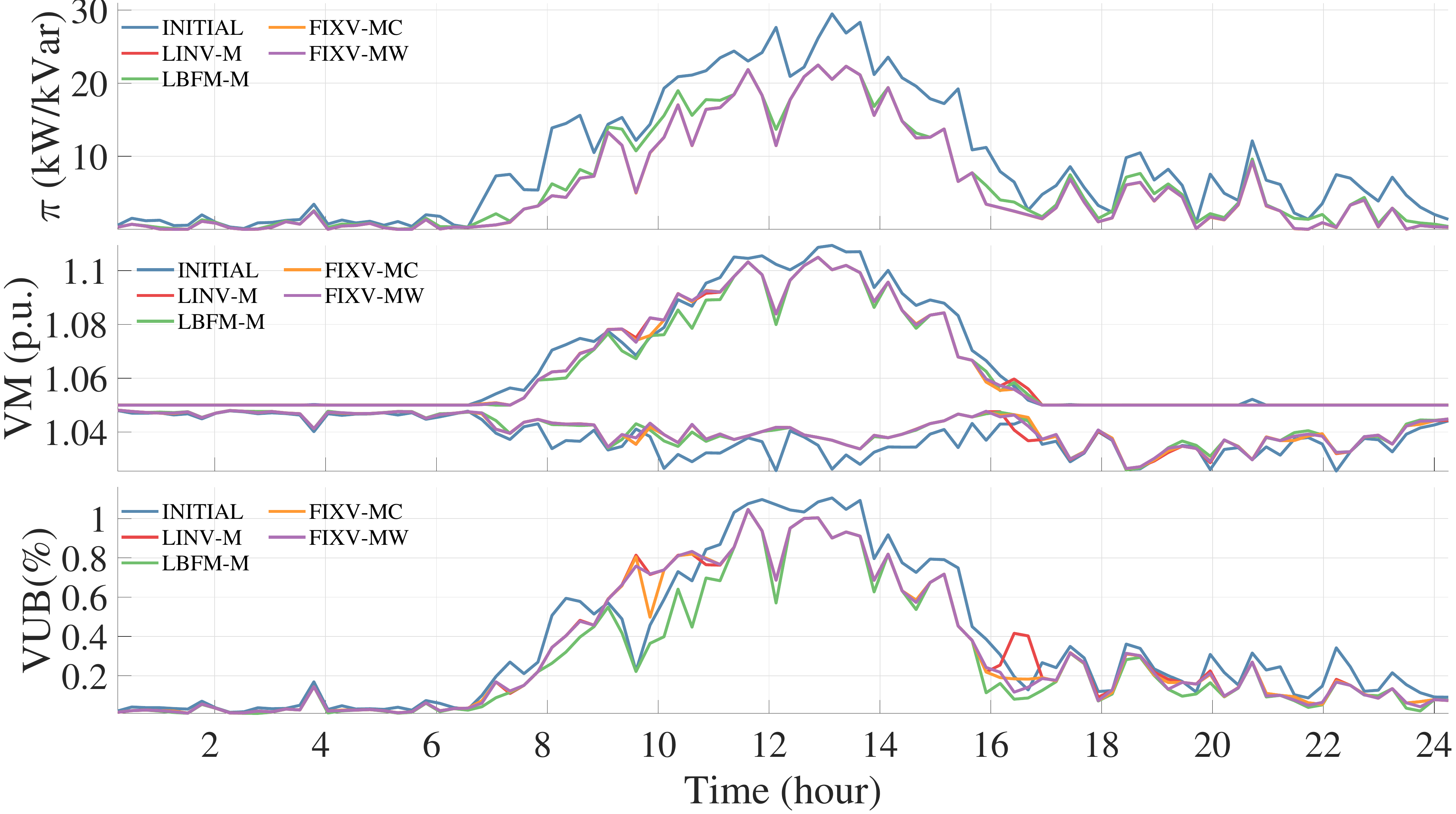}
		\caption{Network Unbalance, VMs range and VUB levels for the Base Case.}
		\label{fig-APPEEC-CaseI-01}
	\end{figure}

	The network unbalances, VMs and VUBs of all periods throughout the day are presented in Fig.\ref{fig-APPEEC-CaseI-01}.
	
	The simulation results clearly show that the system experience high network unbalance, over-voltage issue and VUB violations when PSDs are at original phase connections. By contrast, the network unbalances throughout the day can be effectively mitigated after optimizing PSD phase connections, and the violations of VM and VUB can be effectively addressed or alleviated in some periods. 
	The average network unbalance reduced by 37.27\%, 30.58\%, 37.27\% and 37.19\% under LINV-M, LBFM-M, FIXV-MC and FIXV-MW, respectively. LBFM-M  LINV-M, FIXV-MC and FIXV-MW outperform LBFM-M due to the strong linearizations and assumptions taken in the latter method. It is also noteworthy that although FIXV-MW only iterates once in the algorithm, its performance is close to LINV-M that has been demonstrated accurate and efficient in addressing the network unbalance issue \cite{PRD_LVDN_01,PRD_LVDN_02}.
	
\begin{figure*}[htb!]
		\centering\includegraphics[scale=0.39]{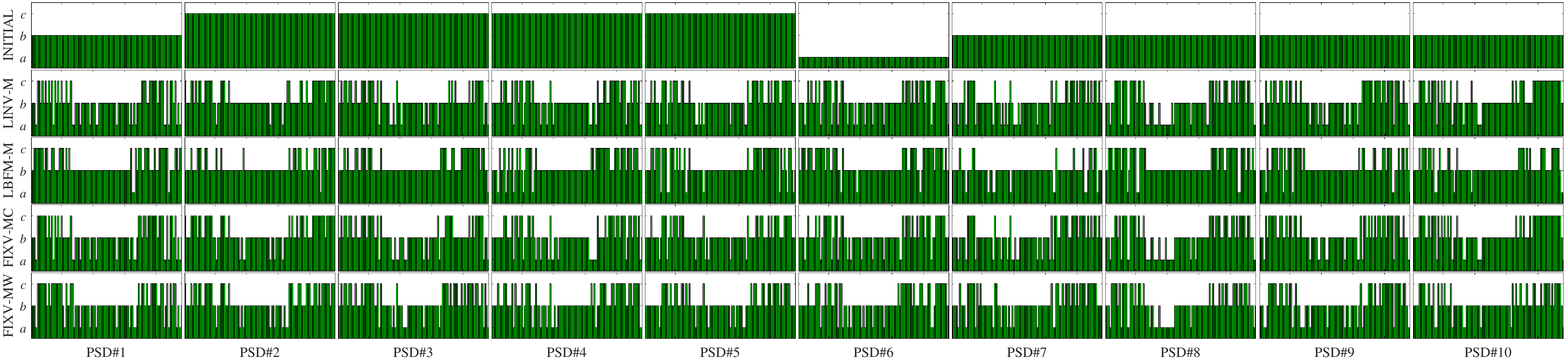}
		\caption{Optimal PSD connections under all approaches.}
		\label{fig-APPEEC-CaseI-04}
	\end{figure*}

    The optimal PSD connections under all approaches are presented in Fig.\ref{fig-APPEEC-CaseI-04}, which shows that the control strategy for LINV-M, FIXV-MC and FIXV-MW are much similar to each other while LBFM-M leads to results with more differences due to the approximations.  

\begin{figure}[htb!]
		\centering\includegraphics[scale=0.23]{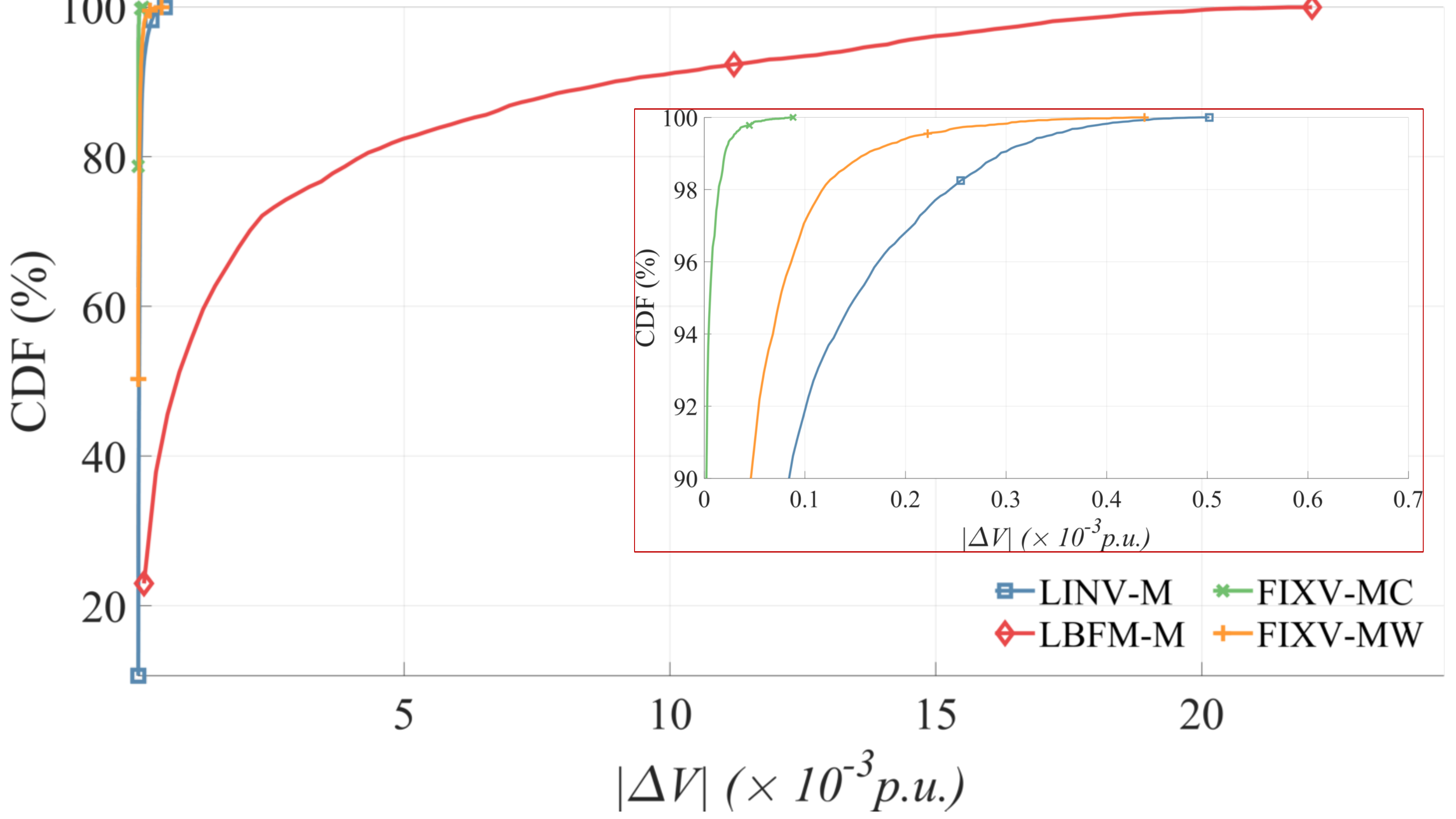}
		\caption{Voltage comparison before/after UTPF for LIVN-M, FIXV-M and LBFM-M.}
		\label{fig-APPEEC-CaseI-02}
	\end{figure}

	To investigate the formulation accuracy, $V_{\phi,i}$ from LINV-M+UTPF, LBFM-M+UTPF, FIXV-MC+UTPF and FIXV-MW+UTPF are compared with their values from LINV-M, LBFM-M, FIXV-MC and FIXV-MW. Specifically, $\Delta V_{\phi,i}$ for each method is presented in Fig.\ref{fig-APPEEC-CaseI-02}. 
	
	$\Delta V_{\phi,i}$ in Fig.\ref{fig-APPEEC-CaseI-02} are presented based on the cumulative distribution function (CDF) under various methods. Any point, say $(x,y)$, on each curve means $100y\%$ of $\Delta V_{\phi,i}$ are less than $x$ under the corresponding method. The CDF figure implies that method $A$ outperforms method $B$ if $A$'s curve is on the left-top side of $B$'s curve. For all the discussed methods, the accuracy of LINV-M is highest, followed by FIXV-MC and then FIXV-MW. As power loss is neglected and angles of nodal voltages are approximated in LBFM-M, the accuracy is much lower than the other methods, and the largest error could reach higher than 0.02 p.u.

	\subsection{The Case with Controllable PV Systems}
	In this case, we assume reactive power of each PV inverter is controllable within [-5\%, 5\%] of its rated capacity. In other words, $Q^n_{i,j}$ will be replaced by a variable $Q^{nx}_{i,j}$ subject to the following constraint if the $j^\text{th}$ customer at node $i$ is with PV. 
    \begin{eqnarray}
    \footnotesize
		\label{QPVcons}
        Q^\text{min}_{i,j}\le Q^{nx}_{i,j}-Q^n_{i,j}\le Q^\text{max}_{i,j}~\forall i,\forall j
	\end{eqnarray}	
	where $Q^\text{max}_{i,j}=-Q^\text{min}_{i,j}=C^{pv}_{i,j}\times 5\%$ with $C^{pv}_{i,j}$ being the PV capacity of the $j^\text{th}$ customer at node $i$.
	
\begin{figure}[htb!]
		\centering\includegraphics[scale=0.2]{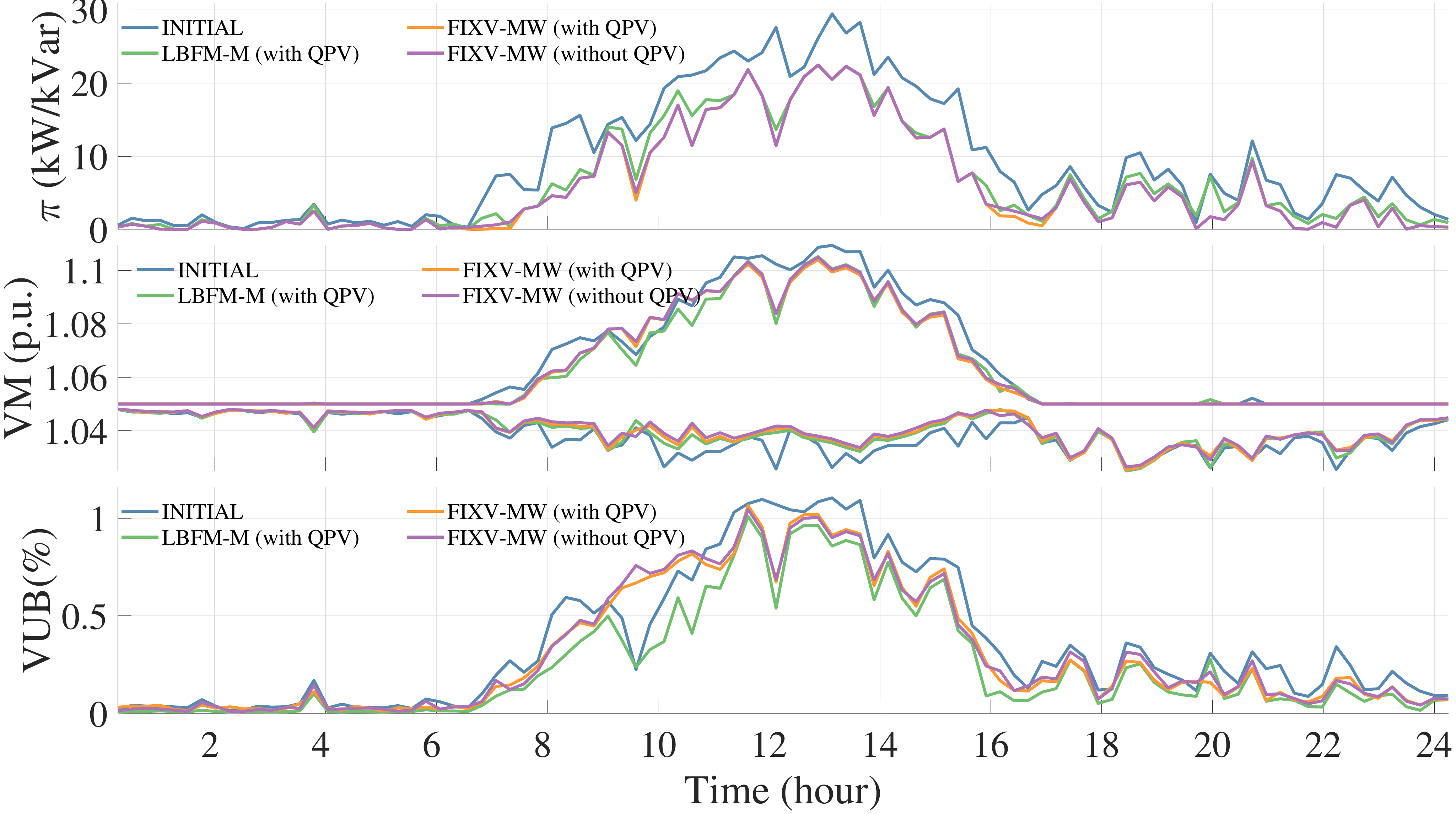}
		\caption{Network Unbalance, VMs range and VUB levels for the Case with controllable PVs (with/without QPV means with/without controllable PVs).}
		\label{fig-APPEEC-CaseII-01}
	\end{figure}

	Average computational time for a single period under FIXV-MW and LBFM-M is 14.22 and 3.39 seconds respectively, demonstrating both of them can flexibly and efficiently deal with more controllable variables. Moreover, the average network unbalance is reduced by 30.25\% and 38.03 for LBFM-M and FIXV-MW, respectively. 
	Simulation results, including network unbalance, VMs and VUBs throughout the day, are presented in Fig.\ref{fig-APPEEC-CaseII-01}.
	
	It should be noted that the simulation results under INITIAL and FIXV-MW from the Base Case are also provided for comparison purposes. Comparing FIXV-MW (without QPV) and FIXV-MW (with QPV) in the figure shows that making PV inverters controllable to provide reactive powers could further help mitigate network unbalance, although the improvement is not significant, which is because reactive power level in the network is much lower than the active power as shown in \ref{fig-system-demand}. Moreover, LBFM-M still leads to a higher network unbalance than IFXV-MW (without QPV), although it utilizes more controllable resources, implying the importance and potential benefits of accurate problem formulation. 
	
	\begin{figure}[htb!]
		\centering\includegraphics[scale=0.23]{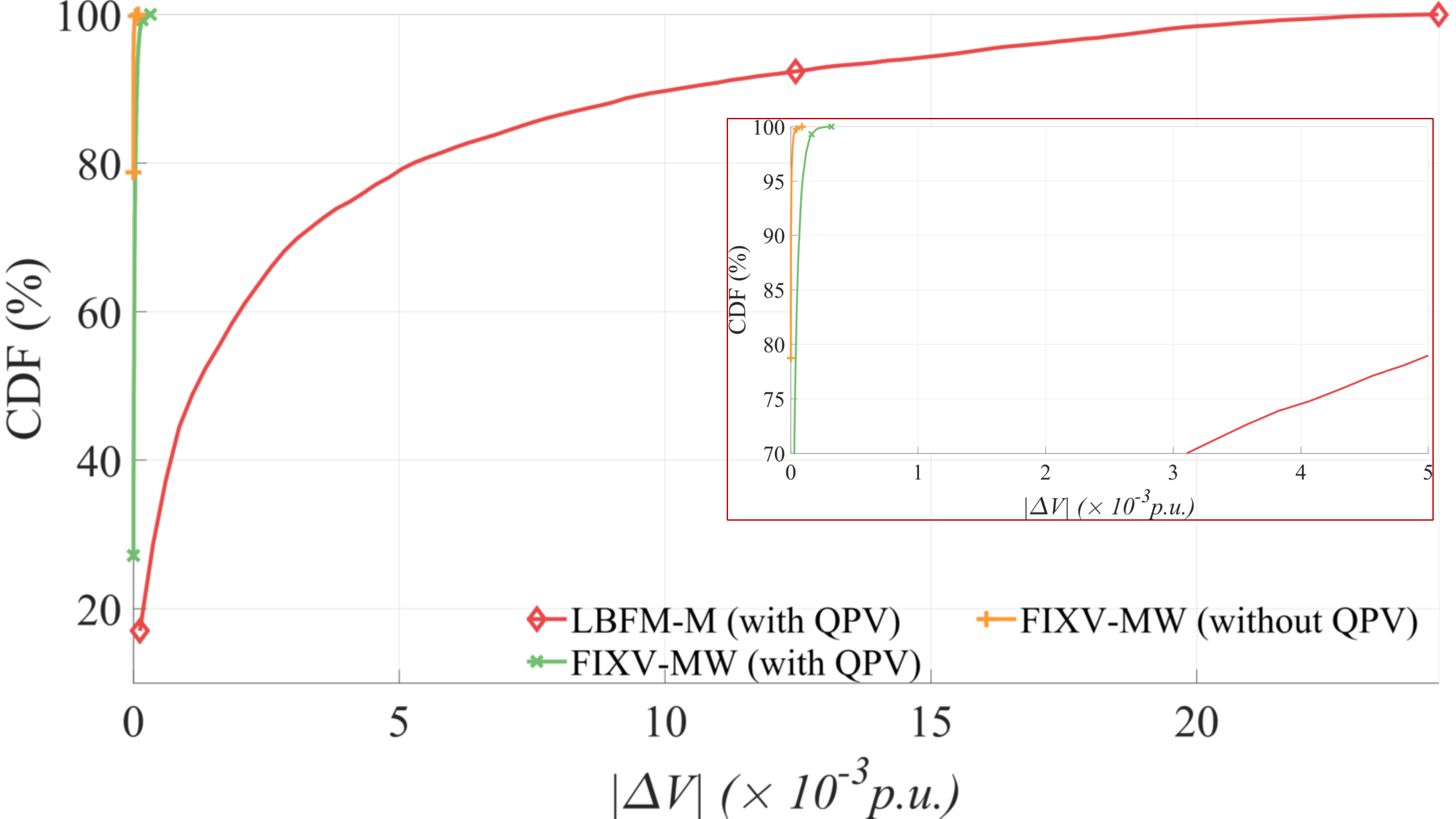}
		\caption{Voltage comparison before/after UTPF for FIXV-MW, FIXV-MW and LBFM-M with/without controllable PVs.}
		\label{fig-APPEEC-CaseII-02}
	\end{figure}

	To further investigate the formulation accuracy, VMs from FIXV-MW+UTPF and LBFM-M+UTPF are compared with VMs directly from FIXV-MW and LBFM-M, which are presented in Fig.\ref{fig-APPEEC-CaseII-02}. Similar to the accuracy analysis in the Base Case, the errors are sufficiently small for FIXV-MW, while the error could be higher than 0.02 p.u. for the LBFM-M method. 
	
	\section{Conclusions}
	In this paper, algorithms for coordinating PSDs to mitigate unbalance in LVDN are studied. Although the previous complete linearization-based method in \cite{PRD_LVDN_01} has been demonstrated effective and accurate, its practicality may be limited due to its computational speed and limited capabilities to deal with more controllable devices. A novel approach based on the well-known linearized UTOPF in \cite{RN64} is demonstrated to be more computationally efficient and flexible. However, the strong assumptions in the linearization limited the accuracy and, as a result, lead to sub-optimal solutions in the real-world. By contrast, the proposed method, which iteratively updates the fixed nodal voltages in addressing the non-convex power balance equations, is shown to be both efficient and flexible. Especially if starting from a warm-start point, the proposed method could achieve higher computational performance while preserving high accuracy, which provides a more efficient approach to coordinate PSDs in real networks. The proposed method can also be adapted in demand response applications, where the underlying three-phase unbalanced network may violate some operational constraints, and hence, become infeasible in practice as discussed in \cite{DR_unbalanced2018}. Optimal demand response by considering operational requirements of an LVDN based on exact formulation of UTOPF and coordination of more controllable devices is the subject of future work. Moreover, as integer variables may be introduced to the UTOPF, developing more efficient second-order cone-based approach to solve UTOPF other than the semi-definite relaxation-based approach as discussed in \cite{RN64}, and investigating its accuracy and computational efficiency also falls in our future research interest.

	\bibliography{OPSDREF_APPEEC}

\begin{thebibliography}{10}
\providecommand{\url}[1]{#1}
\csname url@samestyle\endcsname
\providecommand{\newblock}{\relax}
\providecommand{\bibinfo}[2]{#2}
\providecommand{\BIBentrySTDinterwordspacing}{\spaceskip=0pt\relax}
\providecommand{\BIBentryALTinterwordstretchfactor}{4}
\providecommand{\BIBentryALTinterwordspacing}{\spaceskip=\fontdimen2\font plus
\BIBentryALTinterwordstretchfactor\fontdimen3\font minus
  \fontdimen4\font\relax}
\providecommand{\BIBforeignlanguage}[2]{{%
\expandafter\ifx\csname l@#1\endcsname\relax
\typeout{** WARNING: IEEEtran.bst: No hyphenation pattern has been}%
\typeout{** loaded for the language `#1'. Using the pattern for}%
\typeout{** the default language instead.}%
\else
\language=\csname l@#1\endcsname
\fi
#2}}
\providecommand{\BIBdecl}{\relax}
\BIBdecl

\bibitem{RN61}
P.~K.~C. Wong, A.~Kalam, and R.~Barr, ``Modelling and analysis of practical
  options to improve the hosting capacity of low voltage networks for embedded
  photo-voltaic generation,'' \emph{IET Renew. Power Gener.}, vol.~11, no.~5,
  pp. 625--632, 2017.

\bibitem{RN103}
J.~Zhu, M.-Y. Chow, and F.~Zhang, ``Phase balancing using mixed-integer
  programming,'' \emph{IEEE Trans. Power Syst.}, vol.~13, no.~4, pp.
  1487--1492, 1998.

\bibitem{RN192}
\BIBentryALTinterwordspacing
``Solar report ({January} 2019),'' Australian Energy Council, Report, 2019.
  [Online]. Available:
  \url{https://www.energycouncil.com.au/media/15358/australian-energy-council-solar-report_-january-2019.pdf}
\BIBentrySTDinterwordspacing

\bibitem{BLiu-TUPF}
B.~Liu, K.~Meng, P.~K. Wong, Z.~Y. Dong, C.~Zhang, B.~Wang, T.~Ting, and Q.~Qi,
  ``Improving operation feasibility of low-voltage distribution network by
  phase-switching devices,'' in \emph{Proc. Int. Conf. on Renewable Power
  Generation (IET RPG)}, Shanghai, China, 2019.

\bibitem{psd-knowledge}
\BIBentryALTinterwordspacing
``Jemena {DER} hosting capacity project interim knowledge sharing report,''
  Jemena Electricity Pty. Ltd., Report, 2020. [Online]. Available:
  \url{https://arena.gov.au/assets/2020/08/jemena-der-hosting-capacity-interim-knowledge-sharing-report.pdf}
\BIBentrySTDinterwordspacing

\bibitem{Liu2021}
M.~Z. Liu, L.~N. Ochoa, S.~Riaz, P.~Mancarella, T.~Ting, J.~San, and
  J.~Theunissen, ``{Grid and market services from the edge: using operating
  envelopes to unlock network-aware bottom-Up flexibility},'' \emph{IEEE Power
  and Energy Magazine}, vol.~19, no.~4, pp. 52--62, 2021.

\bibitem{RN71}
F.~Shahnia, P.~J. Wolfs, and A.~Ghosh, ``Voltage unbalance reduction in low
  voltage feeders by dynamic switching of residential customers among three
  phases,'' \emph{IEEE Trans. Smart Grid}, vol.~5, no.~3, pp. 1318--1327, 2014.

\bibitem{RN98}
C.~H. Lin, C.~S. Chen, M.~Y. Huang, H.~J. Chuang, M.~S. Kang, C.~Y. Ho, and
  C.~W. Huang, ``Optimal phase arrangement of distribution feeders using immune
  algorithm,'' in \emph{Proc. Int. Conf. on Intelligent Syst. Applications to
  Power Syst.}, Niigata, Japan, 2007.

\bibitem{RN44}
J.~Horta, D.~Kofman, D.~Menga, and M.~Caujolle, ``Augmenting {DER} hosting
  capacity of distribution grids through local energy markets and dynamic phase
  switching,'' in \emph{Proc. the Ninth Int. Conf. on Future Energy Syst. -
  e-Energy '18}, Karlsruhe, Germany, 2018.

\bibitem{RN105}
X.~Geng, S.~Gupta, and L.~Xie, ``Robust look-ahead three-phase balancing of
  uncertain distribution loads,'' \emph{arXiv preprint arXiv:1810.00425}, 2018.

\bibitem{PRD_LVDN_01}
B.~{Liu}, K.~{Meng}, Z.~Y. {Dong}, P.~K.~C. {Wong}, and T.~{Ting}, ``Unbalance
  mitigation via phase-switching device and static var compensator in
  low-voltage distribution network,'' \emph{IEEE Trans. Power Syst.}, vol.~35,
  no.~6, pp. 4856--4869, 2020.

\bibitem{PRD_LVDN_03}
B.~{Liu}, K.~{Meng}, Z.~{Dong}, P.~K.~C. {Wong}, and W.~{Wei}, ``Optimal
  placement of phase-reconfiguration devices in low-voltage distribution
  network with residential pv generation,'' \emph{IET Renew. Power Gener.}, pp.
  1--1, 2020.

\bibitem{PRD_LVDN_02}
B.~{Liu}, K.~{Meng}, Z.~{Dong}, P.~K.~C. {Wong}, and X.~{Li}, ``Load balancing
  in low-voltage distribution network via phase reconfiguration: An efficient
  sensitivity-based approach,'' \emph{IEEE Trans. Power Del.}, pp. 1--1, 2020.

\bibitem{RN64}
C.~Zhao, E.~Dall'Anese, and S.~H. Low, ``Optimal power flow in multiphase
  radial networks with delta connections,'' in \emph{Proc. IREP 10th Bulk Power
  Syst. Dynamics and Control Symp.}, Espinho, Portugal, 2017.

\bibitem{linear-SOC}
A.~N. Aharon Ben-Tal, ``On polyhedral approximations of the second-order
  cone,'' \emph{Operations Research}, vol.~26, no.~2, pp. 193--205, 2001.

\bibitem{RN38}
H.~Ahmadi, J.~R. Marti, and A.~von Meier, ``A linear power flow formulation for
  three-phase distribution systems,'' \emph{IEEE Trans. Power Syst.}, vol.~31,
  no.~6, pp. 5012--5021, 2016.

\bibitem{RN66}
B.~A. Robbins and A.~D. Dominguez-Garcia, ``Optimal reactive power dispatch for
  voltage regulation in unbalanced distribution systems,'' \emph{IEEE Trans.
  Power Syst.}, vol.~31, no.~4, pp. 2903--2913, 2016.

\bibitem{mccormick}
A.~Gupte, S.~Ahmed, M.~Cheon, and S.~Dey, ``Solving mixed integer bilinear
  problems using {MILP} formulations,'' \emph{{SIAM} Journal on Optimization},
  vol.~23, no.~2, pp. 721--744, 2013.

\bibitem{RN37}
T.~Jen-Hao, ``A direct approach for distribution system load flow solutions,''
  \emph{IEEE Trans. Power Del.}, vol.~18, no.~3, pp. 882--887, 2003.

\bibitem{gurobi}
\BIBentryALTinterwordspacing
L.~Gurobi~Optimization, ``Gurobi optimizer reference manual,'' 2018. [Online].
  Available: \url{http://www.gurobi.com}
\BIBentrySTDinterwordspacing

\bibitem{EuroLVDN}
``Simulation data for the {IEEE} {E}uropean {LVDN},''
  \url{https://site.ieee.org/pes-testfeeders/resources/}.

\bibitem{DR_unbalanced2018}
W.~Zheng, W.~Wu, B.~Zhang, and C.~Lin, ``Distributed optimal residential demand
  response considering operational constraints of unbalanced distribution
  networks,'' \emph{IET Gener. Transm. Distrib.}, vol.~12, no.~9, pp.
  1970--1979, 2018.

\end{thebibliography}
\end{document}